\RequirePackage{ifpdf}
\ifpdf % We are running pdfTeX in pdf mode
\documentclass[pdftex]{sigma}
\else
\documentclass{sigma}
\fi

\numberwithin{equation}{section}
\newtheorem{thm}{Theorem}[section]
\newtheorem{cor}[thm]{Corollary}
\newtheorem{prop}[thm]{Proposition}
\newtheorem{con}[thm]{Conjecture}
\newtheorem{lem}[thm]{Lemma}

\newcommand{\LHS}{\text{LHS}}
\newcommand{\RHS}{\text{RHS}}
\newcommand{\abs}[1]{\lvert#1\rvert}
\newcommand{\la}{\lambda}
\renewcommand{\H}{\mathcal H}
\newcommand{\V}{\mathcal V}
\newcommand{\Z}{\mathbb Z}
\newcommand{\C}{\mathbb C}
\newcommand{\Fe}{\mathbb F}
\newcommand{\Q}{\mathbb Q}
\newcommand{\Sum}[1]{\sum_{\substack{I\subseteq [n] \\ \abs{I}=#1}}}
\newcommand{\nin}{\not\in}
\renewcommand{\L}{\mathcal L}
\newcommand{\R}{\mathcal R}
\newcommand{\F}{\mathcal F}
\newcommand{\Symm}{\mathfrak{S}}

\begin{document}

\allowdisplaybreaks

\renewcommand{\thefootnote}{$\star$}

\renewcommand{\PaperNumber}{055}

\FirstPageHeading

\ShortArticleName{Theta Functions and Kawanaka's Conjecture}

\ArticleName{Theta Functions, Elliptic Hypergeometric Series,\\ and
Kawanaka's Macdonald Polynomial Conjecture\footnote{This paper is a contribution to the Proceedings of the Workshop ``Elliptic Integrable Systems, Isomonodromy Problems, and Hypergeometric Functions'' (July 21--25, 2008, MPIM, Bonn, Germany). The full collection
is available at
\href{http://www.emis.de/journals/SIGMA/Elliptic-Integrable-Systems.html}{http://www.emis.de/journals/SIGMA/Elliptic-Integrable-Systems.html}}}

\Author{Robin LANGER~$^\dag$,
Michael J. SCHLOSSER~$^\ddag$ and S. Ole WARNAAR~$^\S$}

\AuthorNameForHeading{R. Langer, M.J. Schlosser and S.O. Warnaar}

\Address{$^\dag$~Department of Mathematics and Statistics,
The University of Melbourne, VIC 3010, Australia}
\EmailD{\href{mailto:langer.robin@gmail.com}{langer.robin@gmail.com}}

\Address{$^\ddag$~Fakult\"at f\"ur Mathematik, Universit\"at Wien,
Nordbergstrasse 15, A-1090 Vienna, Austria}
\EmailD{\href{mailto:michael.schlosser@univie.ac.at}{michael.schlosser@univie.ac.at}}

\Address{$^\S$~School of Mathematics and Physics, The University of Queensland,\\
\hphantom{$^\S$}~Brisbane, QLD 4072, Australia}
\EmailD{\href{mailto:o.warnaar@maths.uq.edu.au}{o.warnaar@maths.uq.edu.au}}

\ArticleDates{Received March 01, 2009, in f\/inal form May 19,
2009; Published online May 25, 2009}

\Abstract{We give a new theta-function identity,
a special case of which is utilised
to prove Kawanaka's Macdonald polynomial conjecture.
The theta-function identity further yields a~transformation
formula for multivariable elliptic hypergeometric series which
appears to be new even in the one-variable, basic case.}

\Keywords{theta functions; Macdonald polynomials; elliptic hypergeometric series}

\Classification{05E05; 33D52; 33D67}

\renewcommand{\thefootnote}{\arabic{footnote}}
\setcounter{footnote}{0}

\section{Introduction}

The recent discovery of elliptic hypergeometric series (EHS)
by Frenkel and Turaev~\cite{FT97} has led to a renewed interest in
theta-function identities. Such identities are at the core of many of
the proofs of identities for EHS associated with
root systems. If, for $\abs{p}<1$ and $x\in\C^{\ast}$,
\[
\theta(x)=\theta(x;p)=\prod_{k=0}^\infty\big(1-xp^k\big)\big(1-x^{-1}p^{k+1}\big)
\]
is a normalised theta function, then three examples of theta
function identities featured in the theory of EHS are
\begin{gather}\label{WW}
\sum_{i=1}^n\frac{\prod\limits_{j=1}^n\theta(x_i/y_j)}
{\prod\limits_{j=1,\,j\neq i}^n\theta(x_i/x_j)}=0 \qquad\text{for}\quad
x_1\cdots x_n=y_1\cdots y_n,
\\
\label{Gu}
\sum_{i=1}^n\frac{x_i\prod\limits_{j=1}^{n-2}\theta(x_i/y_j)\theta(x_iy_j)}
{\prod\limits_{j=1,\,j\neq i}^n\theta(x_i/x_j)\theta(x_ix_j)}=0
\qquad\text{for}\quad n\geq 2,
\end{gather}
and
\begin{equation}\label{KN}
\Sum{r} \prod_{\substack{i\in I \\ j\in [n]}}
\frac{\theta(x_iy_j)}{\theta(qx_iy_j)}
\prod_{\substack{i\in I \\ j\nin I}}
\frac{\theta(qx_i/x_j)}{\theta(x_i/x_j)}
=\Sum{r} \prod_{\substack{i\in I \\ j\in [n]}}
\frac{\theta(y_ix_j)}{\theta(qy_ix_j)}
\prod_{\substack{i\in I \\ j\nin I}}
\frac{\theta(qy_i/y_j)}{\theta(y_i/y_j)},
\end{equation}
where $[n]:=\{1,2,\dots,n\}$.
The identity \eqref{WW} is standard and can for example be found in
the classic text of Whittaker and Watson \cite[page~451]{WW27}.
It was employed by Gustafson in \cite{Gustafson87} to derive an
$\text{A}_{n-1}$
extension of Bailey's very-well-poised $_6\psi_6$ summation.
In the same paper, Gustafson discovered the identity \eqref{Gu}
\cite[Lemma 4.14]{Gustafson87} and used it to derive a $\text{C}_n$
extension of Bailey's very-well-poised $_6\psi_6$ summation.
The identities \eqref{WW} and \eqref{Gu} were also
employed by Rosengren in \cite{Rosengren04} to prove
elliptic analogues of Milne's $\text{A}_{n-1}$ Jackson sum
\cite{Milne88} and Schlosser's $\text{D}_n$ Jackson sum \cite{Schlosser97}.
The identity \eqref{KN} was only recently
discovered by Kajihara and Noumi \cite[Theorem~1.3]{KN03}
by combining symmetries  of the Frobenius determinant with operator methods.
They used it to extend Rosengren's elliptic $\text{A}_{n-1}$ Jackson summation
to a transformation between EHS on A-type root systems of dif\/ferent rank,
see \cite[Theorem~2.1]{KN03}.
Ruijsenaars, in his study of integrable systems of Calogero--Moser type,
also discovered the identity \eqref{KN}, see
\cite[Equation~(2.44)]{Ruijsenaars05}, where it appeared
in the form of a functional equation encoding the commutativity
of analytic dif\/ference operators of type $\text{A}_{n-1}$.
Earlier, he had used a special case of this identity
(Equation~(2.4) in~\cite{Ruijsenaars87}, see also~\cite{RS86})
to show the commutativity of what is nowadays known as the
Macdonald--Ruijsenaars dif\/ference operator~\cite{Ruijsenaars87}.

Let
\[
(a)_n=(a;q,p)_n=\prod_{k=0}^{n-1}\theta\big(aq^k;p\big)
\]
be a theta shifted factorial (cf.\ \cite[Chapter~11]{GR04}),
and set
\[
(a_1,\dots,a_k)_n=(a_1)_n\cdots (a_k)_n.
\]
Then the main result of this paper is the following new
identity for theta functions.
\begin{thm}\label{Thmrn}
For $n$ a nonnegative integer and $v,w,q,t,x_1,\dots,x_n\in\C^{\ast}$
such that both sides are well defined,
\begin{gather*}
\sum_{r=0}^n \frac{(v,w)_r}{(qv/t,qw/t)_r}
\sum_{\substack{I\subseteq [n] \\ \abs{I}=r}}
\prod_{\substack{i\in I}}
\frac{q \theta(vx_i/tw)\theta(tq^{-r}x_i/qw)}
{t \theta(vx_i/qw)\theta(q^{-r}x_i/w)} \\ \qquad\quad{} \times
\prod_{j\nin I}
\frac{\theta(x_j/q)\theta(q^{-r}x_j/tw)}
{\theta(x_j/t)\theta(q^{-r}x_j/qw)}
\prod_{\substack{i\in I \\ j\nin I}}
\frac{\theta(t x_i/qx_j)\theta(q x_i/x_j)}
{\theta(x_i/x_j)\theta(tx_i/x_j)} \\
\qquad{} =\sum_{r=0}^n \frac{(v,w)_r}{(qv/t,qw/t)_r}
\sum_{\substack{I\subseteq [n] \\ \abs{I}=r}}
\prod_{\substack{i\in I}}
\frac{q \theta(x_i/q)\theta(q^rvx_i/t^2)}
{t \theta(x_i/t)\theta(q^rvx_i/tq)} \\
\qquad\quad{} \times
\prod_{\substack{j\nin I}}
\frac{\theta(vx_j/tw)\theta(q^rvx_j/q)}
{\theta(vx_j/qw)\theta(q^rvx_j/t)}
\prod_{\substack{i\in I \\ j\nin I}}
\frac{\theta(t x_j/q x_i)\theta(q x_j/x_i)}
{\theta(x_j/x_i)\theta(t x_j/x_i)}.
\end{gather*}
\end{thm}
After clearing denominators and replacing
$(t,v,w,x_1)\mapsto (tq,tv,w^{-1},qtx)$
the $n=1$ case of the theorem takes the form
\begin{gather}
\theta(v)\theta(tx)\theta(tw)
\theta(wx)\theta(tvx)\theta\big(t^2vwx\big)
+\theta(w)\theta(x)\theta(tv)
\theta(tvx)\theta\big(t^2wx\big)\theta(tvwx) \nonumber\\
\qquad{}=\theta(v)\theta(x)\theta(tw)
\theta(twx)\theta\big(t^2vx\big)\theta(tvwx)
+\theta(w)\theta(tx)\theta(tv)
\theta(vx)\theta(twx)\theta\big(t^2vwx\big). \!\!\label{nis1}
\end{gather}
This formula, which also follows from the
$(n,r)=(2,1)$ case of \eqref{KN} upon clearing denominators and
substituting $(x_2,y_1,y_2,q)\mapsto (twxx_1,v/x_1,1/txx_1,t)$,
is a once-iterated version of the Riemann
relation\footnote{Riemann gave a number of addition formulae for
theta functions on arbitrary genus Riemann surfaces, but it is
not clear he actually was the f\/irst to discover \eqref{rr}.}
\begin{gather}\label{rr}
\theta(xz)\theta(x/z)\theta(yw)\theta(y/w)-
\theta(xw)\theta(x/w)\theta(yz)\theta(y/z)
=\frac{y}{z} \theta(xy)\theta(x/y)\theta(zw)\theta(z/w).
\end{gather}

We remark that formula \eqref{nis1} is a four-term identity involving
\textit{four} free parameters, each term containing a product of
\textit{six} theta functions.
For comparison, the formula \eqref{WW} (\eqref{Gu}) for $n=4$
is a four-term identity
involving \textit{seven} (\textit{six}) free parameters,
each term containing a~product of \textit{seven} (\textit{ten})
theta functions.

As will be shown in Section~\ref{secell}, Theorem~\ref{Thmrn}
can be used to obtain identities for EHS for the root systems of type~A.
More unexpectedly, however, it also implies a combinatorial identity
in the theory of Macdonald polynomials \cite[Chapter~VI]{Macdonald95},
conjectured in 1999 by Kawanaka~\cite{Kawanaka99}.
For~$\lambda$ a partition let $P_{\la}(x;q,t)$ be the Macdonald
symmetric function in countably many independent variables
$x=(x_1,x_2,\dots)$.
Let $a(s)$ and $l(s)$ be the arm-length and leg-length of the
square $s$ in the diagram of $\la$, and let
\[
(a;q)_{\infty}=\prod_{k=0}^{\infty}\big(1-aq^k\big)
\]
be a $q$-shifted factorial.
\begin{con}[Kawanaka]\label{con}
The following formal identity holds:
\begin{equation}\label{Kawanaka}
\sum_{\la} \prod_{s\in\la}
\biggl(\frac{1+q^{a(s)}t^{l(s)+1}}{1-q^{a(s)+1}t^{l(s)}}\biggr)
P_{\la}\big(x;q^2,t^2\big)
=\prod_{i\geq 1}\frac{(-tx_i;q)_{\infty}}{(x_i;q)_{\infty}}
\prod_{i<j}\frac{(t^2x_ix_j;q^2)_{\infty}}
{(x_ix_j;q^2)_{\infty}}.
\end{equation}
\end{con}
Using elementary results from Macdonald polynomial theory
and (a special limiting case of) Theorem~\ref{Thmrn}
we can claim the following.
\begin{thm}\label{thmcon}
Kawanaka's conjecture is true.
\end{thm}

We f\/inally remark that although Theorem~\ref{Thmrn} is new, a special
limiting case coincides with a~limiting case of another theta-function
identity, implicit in \cite[Corollary 5.3]{Rosengren04}.
\begin{thm}[Rosengren]\label{ThmR}
For $n$ a nonnegative integer and $v,w,y,z,q,x_1,\dots,x_n\in \C^{\ast}$
such that $vw=q^{n-1}yz$ and such that both sides are well defined,
\begin{gather*}
\sum_{r=0}^n (-1)^r q^{\binom{r+1}{2}-nr}
\frac{(v,w)_r}{(y,z)_r}
\sum_{\substack{I\subseteq [n] \\ \abs{I}=r}}
\prod_{i\in I} \frac{\theta(yx_i)\theta(zx_i)}
{\theta(q^{1-r}x_i)}
\prod_{j\nin I}\frac{\theta(vx_j)\theta(wx_j)}{\theta(q^{-r}x_j)}
\prod_{\substack{i\in I \\ j\nin I}}
\frac{\theta(qx_i/x_j)}{\theta(x_i/x_j)}  \\
\qquad{} =\frac{(y/v,y/w)_n}{(y,y/vw)_n}
\prod_{i=1}^n \theta(vwx_i).
\end{gather*}
\end{thm}
If we set $p=0$ in Theorems~\ref{Thmrn} and Theorem~\ref{ThmR}
using $\theta(a;0)=(1-a)$, and then let $t\to\infty$ in the former
(so that its right-hand side vanishes unless $r=0$) and
$(y,z)\to (0,\infty)$ (such that $yz=q^{1-n}vw$) in the latter,
we obtain one and the same rational function identity
(up to a~rescaling of $x_i\mapsto wqx_i$ in Theorem~\ref{Thmrn}).

\section{Proof of Theorem~\ref{Thmrn}}

\subsection{Preliminary remarks}
Let the left- and right-hand sides of the identity
of the theorem be denoted as $L(x;v,w,q,t;p)$ and
$R(x;v,w,q,t;p)$ respectively, where $x:=(x_1,\dots,x_n)$.
Then, by
\begin{equation}\label{thetainv}
\theta(z)=-z \theta\big(z^{-1}\big)
\end{equation}
and
\[
\frac{(a;q^{-1},p)_k}{(b;q^{-1},p)_k}=
\frac{(a^{-1};q,p)_k}{(b^{-1};q,p)_k}\left(\frac{a}{b}\right)^k,
\]
it follows that
\begin{equation}\label{symm}
L(x;v,w,q,t;p)=
R\big(q^{-1}t^{-1}vw^{-1}x;v^{-1},w^{-1},q^{-1},t^{-1};p\big).
\end{equation}
Hence Theorem~\ref{Thmrn} may be reformulated as the symmetry
\[
L(x;v,w,q,t;p)=L\big(q^{-1}t^{-1}vw^{-1}x;v^{-1},w^{-1},q^{-1},t^{-1};p\big).
\]
Alternatively, by
\[
\frac{(a)_{n-k}}{(b)_{n-k}}=
\frac{(a)_n(q^{1-n}/b)_k}{(b)_n(q^{1-n}/a)_k}
\left(\frac{b}{a}\right)^k,
\]
and the substitutions $I\mapsto [n]-I$ and $r\mapsto n-r$, it follows that
\[
L(x;v,w,q,t;p)=R\big(x;q^{-n}tw^{-1},q^{-n}tv^{-1},q,t;p\big)
\frac{(v,w)_n}{(qv/t,qw/t)_n} \left(\frac{q}{t}\right)^n.
\]
Hence Theorem~\ref{Thmrn} is also equivalent to
\[
L(x;v,w,q,t;p)=L\big(x;q^{-n}tw^{-1},q^{-n}tv^{-1},q,t;p\big)
\frac{(v,w)_n}{(qv/t,qw/t)_n} \left(\frac{q}{t}\right)^n.
\]

\subsection{Proof of Theorem~\ref{Thmrn}}
Recall that $x=(x_1,\dots,x_n)$.
We begin by introducing a scalar variable $u$ in the
theorem by making the substitution $x\mapsto x/u$. Then
\begin{gather*}
\sum_{r=0}^n \frac{(v,w)_r}{(qv/t,qw/t)_r}
\sum_{\substack{I\subseteq [n] \\ \abs{I}=r}}
\prod_{\substack{i\in I}}
\frac{q \theta(vx_i/tuw)\theta(tq^{-r}x_i/quw)}
{t \theta(vx_i/quw)\theta(q^{-r}x_i/uw)} \\ \qquad\quad{} \times
\prod_{j\nin I}
\frac{\theta(x_j/qu)\theta(q^{-r}x_j/tuw)}
{\theta(x_j/tu)\theta(q^{-r}x_j/quw)}
\prod_{\substack{i\in I \\ j\nin I}}
\frac{\theta(t x_i/qx_j)\theta(q x_i/x_j)}
{\theta(x_i/x_j)\theta(tx_i/x_j)} \\
\qquad{}=\sum_{r=0}^n \frac{(v,w)_r}{(qv/t,qw/t)_r}
\sum_{\substack{I\subseteq [n] \\ \abs{I}=r}}
\prod_{\substack{i\in I}}
\frac{q \theta(x_i/qu)\theta(q^rvx_i/t^2u)}
{t \theta(x_i/tu)\theta(q^rvx_i/tqu)} \\
\qquad\quad{} \times
\prod_{\substack{j\nin I}}
\frac{\theta(vx_j/tuw)\theta(q^rvx_j/qu)}
{\theta(vx_j/quw)\theta(q^rvx_j/tu)}
\prod_{\substack{i\in I \\ j\nin I}}
\frac{\theta(t x_j/q x_i)\theta(q x_j/x_i)}
{\theta(x_j/x_i)\theta(t x_j/x_i)}.
\end{gather*}
Let $\L(x;u,v,w,q,t;p)$ and $\R(x;u,v,w,q,t;p)$ denote the
left-hand side and right-hand side of this identity,
and further def\/ine
\[
\F=\L-\R.
\]
Comparing with our earlier def\/initions we thus have
\begin{gather*}
\L(x;u,v,w,q,t;p)=L(x/u;v,w,q,t;p) \quad \ \text{and} \ \quad
\R(x;u,v,w,q,t;p)=R(x/u;v,w,q,t;p).\!
\end{gather*}
We are mainly interested in the $u$-dependence of $\L$, $\R$ and
$\F$, and will frequently write $\L(u)$, $\R(u)$ and
$\F(u)$.
The claim of the theorem is thus $\F(u)=0$, which will be
proved by induction on $n$, the cardinality of the alphabet $x$.

For $n=0$ the theorem is trivial: $\F(\,\text{--}\,;u,v,w,q,t;p)=1-1=0$.

\medskip

From $\theta(z)=-z \theta(pz)$
it immediately follows that $\F$ is periodic along
annuli, with period $p$:
\begin{equation}\label{periodicity}
\F(pu)=\F(u).
\end{equation}
The function $\F(u)$ has simple poles at
\begin{equation}\label{poles}
u=\begin{cases}
x_ip^k/t, &  \\[1mm]
vx_ip^k/qw, &  \\[1mm]
q^{-m}p^k x_i/w & \text{for $m\in[r]$}, \\[1mm]
q^{m-1}p^k v x_i/t & \text{for $m\in[r]$},
\end{cases}
\end{equation}
where $i\in[n]$ and $k\in\Z$.
If we show that these poles all have zero residue then,
by \eqref{periodicity} and Liouville's theorem,
$\F(u)$ must be constant.
By the periodicity \eqref{periodicity} it suf\/f\/ices to
consider the poles \eqref{poles} with $k=0$.
By the symmetry of $\F(x;u,v,w,q,t;p)$ in $x$
it also suf\/f\/ices to only consider $i=n$. From
\eqref{symm} it follows that
\[
\L(x;u,v,w,q,t;p)=\R\big(x;qtuv^{-1}w,v^{-1},w^{-1},q^{-1},t^{-1};p\big).
\]
It is thus suf\/f\/icient to only consider the residues of the poles at
\[
u=\begin{cases}
x_n/t, &  \\[1mm]
q^{-m}x_n/w & \text{for $m\in[r]$.}
\end{cases}
\]

First we compute the residue at $u=x_n/t$. For $\L(u)$
only terms such that $n\nin I$ contribute
and for $\R(u)$ only terms with $n\in I$ do.
Hence, after an elementary calculation
which includes the use of \eqref{thetainv},
\[
\lim_{u\to a}\frac{u-a}{\theta(a/u;p)}=\frac{a}{(p;p)_{\infty}^2}
\]
and a shift $r\mapsto r+1$ in $\R(u)$, we f\/ind
\begin{gather*}
\mathop{\rm Res}_{u=t^{-1}x_n}\F(x;u,v,w,q,t;p) \\
\qquad{}=\F\big(x^{(n)};x_n/qt,v,qw,q,t;p\big)
\frac{qx_n}{t^2} \frac{1}{(p;p)_{\infty}^2}
\frac{\theta(w)\theta(t/q)}{\theta(qw/t)}
\prod_{i=1}^{n-1}
\frac{\theta(q x_i/x_n)\theta(t x_i/qx_n)}
{\theta(tx_i/x_n)\theta(x_i/x_n)},
\end{gather*}
where $x^{(n)}:=(x_1,\dots,x_{n-1})$.
By induction on $n$ this vanishes.

Next we consider the pole at $u=q^{-m}x_n/w$.
The only contributions to its residue
come from~$\L(u)$ with (i) $n\in I$ and $r=m$ or (ii)
$n\nin I$ and $r=m-1$. An elementary calculation shows that
these two contributions are the same up to a sign, and thus
cancel.

Now that we have established that all poles of $\F(u)$
have zero residue we may conclude that $\F(u)$ is independent of $u$.
To show that it is actually identically zero we take $u=x_n/q$.
In $\L(x_n/q)$ only terms such that $n\in I$ contribute
and in $\R(x_n/q)$ only terms with $n\nin I$ do.
Again using \eqref{thetainv}
and making a shift $r\mapsto r+1$ in $\L(x_n/q)$, we f\/ind
\[
\F(x;x_n/q,v,w,q,t;p)=\F\big(x^{(n)};x_n,qv,w,q,t;p\big)
\frac{\theta(v)\theta(qv/tw)}{\theta(qv/t)\theta(v/w)}.
\]
By induction this once again vanishes.\hfill\qed

\section{Proof of Kawanaka's conjecture}

\subsection{Preliminary remarks}
Kawanaka's identity complements a set of four Macdonald
polynomial identities discovered by Macdonald
\cite[page~349]{Macdonald95}. In slightly more general form
as given in \cite{Warnaar06}, these identities may be stated as
the following pair of results (Macdonald's formulae correspond
to $b=0$ and $b=1$)
\begin{subequations}\label{MD}
\begin{gather}
\sum_{\la} b^{c(\la)} \!\!\prod_{\substack{s\in\la \\ l(s)\text{ even}}}\!\!
 \biggl(\frac{1-q^{a(s)}t^{l(s)+1}}{1-q^{a(s)+1}t^{l(s)}}\biggr)
P_{\la}(x;q,t)
=\prod_{i\geq 1}\!\frac{(btx_i;q)_{\infty}}{(bx_i;q)_{\infty}}
\prod_{i<j}\!\frac{(tx_ix_j;q)_{\infty}}{(x_ix_j;q)_{\infty}},\!\! \\
\sum_{\la} b^{r(\la)} \!\! \prod_{\substack{s\in\la \\ a(s)\text{ odd}}}\!\!
 \biggl(\frac{1-q^{a(s)}t^{l(s)+1}}{1-q^{a(s)+1}t^{l(s)}}\biggr)
P_{\la}(x;q,t)
=\prod_{i\geq 1}\!\frac{(1+bx_i)(qtx_i^2;q^2)_{\infty}}
{(x_i^2;q^2)_{\infty}}
\prod_{i<j}\!\frac{(tx_ix_j;q)_{\infty}}{(x_ix_j;q)_{\infty}},\!\!
\end{gather}
\end{subequations}
where $c(\la)$ and $r(\la)$ are the number of columns and rows of odd length,
respectively.
Due to its quadratic nature Kawanaka's identity is signif\/icantly
harder to prove than \eqref{MD}.

If $x$ contains a single variable then \eqref{Kawanaka} simplif\/ies to
the classical $q$-binomial theorem \cite[Equation~(II.3)]{GR04}
\begin{equation}\label{qbt}
\sum_{k=0}^{\infty}\frac{(-t;q)_k}{(q;q)_k}  x^k =
\frac{(-tx;q)_{\infty}}{(x;q)_{\infty}},
\end{equation}
where, for integer $k$,
\[
(a;q)_k=\frac{(a;q)_\infty}{(a q^{k};q)_\infty}.
\]

If $q=t$ then \eqref{Kawanaka} reduces to an identity for
the Schur function $s_{\la}$ proved by Kawanaka \cite{Kawanaka99}.
Specif\/ically, using that $P_{\la}(x;q,q)=s_{\la}(x)$
and $a(s)+l(s)+1=h(s)$ with $h(s)$ the hook-length of the
square $s$, it follows that the $q=t$ case of \eqref{Kawanaka} is
\begin{equation}\label{Kqt}
\sum_{\la} \prod_{s\in\la}
\biggl(\frac{1+q^{h(s)}}{1-q^{h(s)}}\biggr) s_{\la}(x)
=\prod_{i\geq 1}\frac{(-qx_i;q)_{\infty}}{(x_i;q)_{\infty}}
\prod_{i<j}\frac{1}{1-x_ix_j}.
\end{equation}
This result was reproved and reinterpreted by Rosengren in \cite{Rosengren08}.
If $Q_{\mu}$ is Schur's $Q$-function, see \cite[Section~III.8]{Macdonald95},
then Rosengren observed that for $\mu$ a partition of length $m$
\begin{gather*}
Q_{\mu}(1,q,q^2,\dots)
=\biggl(\frac{(-1;q)_{\infty}}{(q;q)_{\infty}}\biggr)^m\!\!\!
\prod_{1\leq i<j\leq m}\!\!(q^{\mu_i}-q^{\mu_j})
\sum_{\la} \prod_{s\in\la}
\biggl(\frac{1+q^{h(s)}}{1-q^{h(s)}}\biggr)
s_{\la}\bigl(-q^{\mu_1},\dots,-q^{\mu_m}\bigr).\!\!
\end{gather*}
By \eqref{Kqt} this results in a product-form for
$Q_{\mu}(1,q,q^2,\dots)$ and, consequently, in a
product-form for the generating function of marked shifted tableaux
\cite[Corollary~3.1]{Rosengren08}. It is an open problem to f\/ind a
corresponding interpretation of~\eqref{Kawanaka}.

Another special case of \eqref{Kawanaka} proved by Kawanaka
corresponds to $q=0$ \cite{Kawanaka91}.
Using $P_{\la}(x;0,t)=P_{\la}(x;t)$ with on the right
a Hall--Littlewood symmetric function, it follows that the
$q=0$ case of \eqref{Kawanaka} is
\[
\sum_{\la} \biggl(\prod_{i\geq 1}\big(-t;t^2\big)_{m_i(\la)}\biggr)
P_{\la}\big(x;t^2\big)
=\prod_{i\geq 1}\frac{1+tx_i}{1-x_i}
\prod_{i<j}\frac{1-t^2x_ix_j}{1-x_ix_j},
\]
with $m_i(\la)$ the multiplicity of the part $i$ in $\lambda$.
This is in fact a special case of a much more general identity
for Hall--Littlewood functions proved in \cite[Theorem 1.1]{Warnaar06}.
So far our attempts to generalise \eqref{Kawanaka} to include this
more general result for Hall--Littlewood functions have been unsuccessful.

\subsection{Macdonald polynomials}

Let $\la=(\la_1,\la_2,\dots)$ be a partition,
i.e., $\la_1\geq \la_2\geq \cdots$ with f\/initely many $\la_i$
unequal to zero.
The length and weight of~$\la$, denoted by
$l(\la)$ and $\abs{\la}$, are the number and sum
of the non-zero~$\la_i$ respectively.
As usual we identify two partitions that dif\/fer only in their
string of zeros, so that $(6,3,3,1,0,0)$ and $(6,3,3,1)$ represent the
same partition.
When $\abs{\la}=N$ we say that~$\la$ is a~partition of $N$,
and the unique partition of zero is denoted by~$0$.
The multiplicity of the part $i$ in the partition~$\la$ is denoted
by $m_i(\la)$.

We identify a partition with its Ferrers graph,
def\/ined by the set of points in $(i,j)\in \Z^2$ such that
$1\leq j\leq \la_i$, and further make the usual
identif\/ication between Ferrers graphs and (Young) diagrams
by replacing points by squares.
The conjugate $\la'$ of $\la$ is the partition obtained by
ref\/lecting the diagram of $\la$ in the main diagonal.

The dominance partial order on the set of partitions of $N$ is
def\/ined by $\la\geq \mu$ if
$\la_1+\cdots+\la_i\geq \mu_1+\cdots+\mu_i$ for all $i\geq 1$.
If $\la\geq \mu$ and $\la\neq\mu$ then $\la>\mu$.

If $\la$ and $\mu$ are partitions then $\mu\subseteq\la$
if (the diagram of) $\mu$ is contained in (the diagram of)~$\la$, i.e., $\mu_i\leq\la_i$ for all $i\geq 1$.
If $\mu\subseteq\la$ then the skew-diagram $\la-\mu$
denotes the set-theoretic dif\/ference between $\la$ and $\mu$,
i.e., those squares of $\la$ not contained in $\mu$.
The skew diagram $\la-\mu$ is a vertical $r$-strip if
$\abs{\la-\mu}:=\abs{\la}-\abs{\mu}=r$ and if, for
all $i\geq 1$, $\la_i-\mu_i$
is at most one, i.e., each row of $\la-\mu$ contains at most
one square. For example, if $\la=(5,4,2,2,1)$ and $\mu=(4,3,1,1,1)$
then $\la-\mu$ is a vertical $4$-strip.
The set of all vertical $r$-strips is denoted by
$\V_r$ and the set of all vertical strips by
$\V=\bigcup_{r=0}^{\infty}\V_r$.
The skew diagram $\la-\mu$ is a horizontal $r$-strip if
$\abs{\la-\mu}=r$ and if,
for all $i\geq 1$, $\la'_i-\mu'_i$
is at most one, i.e., each column of $\la-\mu$ contains at most
one square. The set of all horizontal $r$-strips is denoted
by $\H_r$ and the set of all horizontal strips by $\H$.

Let $s=(i,j)$ be a square in the diagram of $\la$, and let
$a(s)$ and $l(s)$ be the arm-length and leg-length of $s$,
given by
\[
a(s)=\la_i-j, \qquad l(s)=\la'_j-i.
\]
Then we def\/ine the rational functions $b^{+}_{\la}(q,t)$ and
$b^{-}_{\la}(q,t)$ as
\begin{equation}\label{bdef}
b^{\pm}_{\la}(q,t)=\prod_{s\in\la}\frac{1 \mp q^{a(s)}t^{l(s)+1}}
{1-q^{a(s)+1}t^{l(s)}}.
\end{equation}
The function $b^{+}_{\la}(q,t)$ is standard in Macdonald polynomial
theory and is usually denoted as $b_{\la}(q,t)$. Below we use both
notations: $b_{\la}=b_{\la}^{+}$.
Note that $b^{-}_{\la}(q,t)$ corresponds to the product
in the summand of Kawanaka's conjecture.
Since under conjugation arms and legs are interchanged,
it easily follows that
\begin{equation}\label{bp}
b^{-}_{\la'}(q,t)=b^{-}_{\la}(t,q)/b_{\la}\big(t^2,q^2\big).
\end{equation}

Subsequently we require non-combinatorial expressions for both
$b^{-}_{\la}$ and $b^{-}_{\la'}$.
{}From \eqref{bdef} it follows that
\begin{equation}\label{Iexp}
b^{\pm}_{\la}(q,t)=\prod_{i=1}^n
\frac{(\pm t^{n-i+1};q)_{\la_i}}{(qt^{n-i};q)_{\la_i}}
\prod_{1\leq i<j\leq n}
\frac{(\pm t^{j-i},qt^{j-i};q)_{\la_i-\la_j}}
{(\pm t^{j-i+1},qt^{j-i-1};q)_{\la_i-\la_j}},
\end{equation}
where $n$ is an integer such that $l(\la)\le n$. From this and
\eqref{bp} we also f\/ind
\begin{equation}\label{bprimeexp}
b^{-}_{\la'}(q,t)=
\prod_{i=1}^n \frac{(-tq^{n-i};t)_{\la_i}}{(q^{n-i+1};t)_{\la_i}}
\prod_{1\leq i<j\leq n}
\frac{(q^{j-i+1},-tq^{j-i-1};t)_{\la_i-\la_j}}
{(q^{j-i},-tq^{j-i};t)_{\la_i-\la_j}},
\end{equation}
where, again, $l(\la)\le n$.

Let $\Symm_n$ denote the symmetric group, acting on
$x=(x_1,\dots,x_n)$ by permuting the $x_i$,
and let $\Lambda_n=\Z[x_1,\dots,x_n]^{\Symm_n}$ and
$\Lambda$ denote the ring of symmetric polynomials in $n$ independent
variables and the ring of symmetric functions in countably many variables,
respectively.

For $\la=(\la_1,\dots,\la_n)$ a partition
of at most $n$ parts the monomial symmetric function $m_{\la}$
is def\/ined as
\begin{equation*}
m_{\la}(x)=\sum x^{\alpha},
\end{equation*}
where the sum is over all distinct permutations $\alpha$ of
$\la$, and $x^{\alpha}=x_1^{\alpha_1}\cdots x_n^{\alpha_n}$.
For $l(\la)>n$ we set $m_{\la}(x)=0$.
The monomial symmetric functions $m_{\la}$ for $l(\la)\leq n$
form a $\Z$-basis of $\Lambda_n$.

For $r$ a nonnegative integer the power sums $p_r$ are given by
$p_0=1$ and $p_r=m_{(r)}$ for $r>1$. Hence
\begin{equation*}%\label{powersums}
p_r(x)=\sum_{i\geq 1} x_i^r.
\end{equation*}
More generally the power-sum products are def\/ined as
$p_{\la}(x)=p_{\la_1}(x)\cdots p_{\la_n}(x)$.

Def\/ine the Macdonald scalar product $\langle \cdot,\cdot \rangle_{q,t}$
on the ring of symmetric functions by
\begin{equation*}
\langle p_{\la},p_{\mu}\rangle_{q,t}=
\delta_{\la\mu} z_{\la} \prod_{i=1}^n
\frac{1-q^{\la_i}}{1-t^{\la_i}},
\end{equation*}
with $z_{\la}=\prod_{i\geq 1} m_i! \: i^{m_i}$ and
$m_i=m_i(\la)$. If we denote the ring of symmetric functions
in $n$ variables over the f\/ield $\Fe=\Q(q,t)$ of rational functions
in $q$ and $t$ by $\Lambda_{n,\Fe}$, then
the Macdonald polynomial $P_{\la}(x;q,t)$
is the unique symmetric polynomial in $\Lambda_{n,\Fe}$ such that
\cite[Section~VI.4, Equation~(4.7)]{Macdonald95}:
\begin{equation*}%\label{Pm}
P_{\la}(x;q,t)=m_{\la}(x)+\sum_{\mu<\la}
u_{\la\mu}(q,t) m_{\mu}(x)
\end{equation*}
and
\begin{equation*}
\langle P_{\la},P_{\mu} \rangle_{q,t}
=0\quad \text{if$\quad\la\neq\mu$.}
\end{equation*}
The Macdonald polynomials $P_{\la}(x;q,t)$ with $l(\la)\leq n$
form an $\Fe$-basis of $\Lambda_{n,\Fe}$. If $l(\la)>n$ then
$P_{\la}(x;q,t)=0$.

Since $P_{\la}(x_1,\dots,x_n,0;q,t)=P_{\la}(x_1,\dots,x_n;q,t)$
one can extend the Macdonald polynomials to symmetric functions
containing an inf\/inite number of independent variables $x=(x_1,x_2,\dots)$,
to obtain a basis of $\Lambda_{\Fe}=\Lambda\otimes\Fe$.

A second Macdonald symmetric function is def\/ined as
\[
Q_{\la}(x;q,t)=b_{\la}(q,t)P_{\la}(x;q,t).
\]
The normalisation of the Macdonald inner product is then
\begin{equation*}
\langle P_{\la},Q_{\la}\rangle_{q,t}=1.
\end{equation*}

Important in the proof of Kawanaka's conjecture are the Pieri
rules for Macdonald polynomials.
Let $g_r(x;q,t):=Q_{(r)}(x;q,t)$, or equivalently,
\cite[Section~VI.2, Equation~(2.8)]{Macdonald95}
\begin{equation}\label{grprod}
\prod_{i=1}^n\frac{(tx_iy;q)_{\infty}}{(x_iy;q)_{\infty}}
=\sum_{r=0}^{\infty} g_r(x;q,t) y^r.
\end{equation}
Then the Pieri coef\/f\/icients $\phi_{\mu/\nu}$ and $\psi_{\mu/\nu}$
are given by \cite[Section~VI.6, Equation~(6.24)]{Macdonald95}
\begin{subequations}
\begin{gather}\label{P1}
P_{\nu}(x;q,t) g_r(x;q,t) =
\sum_{\substack{\mu \\ \mu-\nu\in \H_r}} \phi_{\mu/\nu}(q,t) P_{\mu}(x;q,t),\\
\label{P2}
Q_{\nu}(x;q,t) g_r(x;q,t) =
\sum_{\substack{\mu \\ \mu-\nu\in \H_r}} \psi_{\mu/\nu}(q,t) Q_{\mu}(x;q,t).
\end{gather}
\end{subequations}

\subsection{Proof of Theorem~\ref{thmcon}}
By \eqref{bdef} Kawanaka's conjecture can be stated as
\begin{equation}\label{idI}
\sum_{\la} b^{-}_{\la}(q,t) P_{\la}\big(x;q^2,t^2\big)
=\prod_{i\geq 1}\frac{(-tx_i;q)_{\infty}}{(x_i;q)_{\infty}}
\prod_{i<j}\frac{\big(t^2x_ix_j;q^2\big)_{\infty}}
{(x_ix_j;q^2)_{\infty}}.
\end{equation}

Our initial steps closely follow Macdonald's proof of \eqref{MD}.
It suf\/f\/ices to prove \eqref{idI} for the f\/inite
set of variables $x=(x_1,\dots,x_n)$.
If we also denote $x'=(x_1,\dots,x_n,y)$ and let
\begin{equation}\label{Phi}
\Phi(x;q,t):=\sum_{\la} b^{-}_{\la}(q,t) P_{\la}\big(x;q^2,t^2\big)
\end{equation}
then, by induction on $n$, it is enough to prove that
\begin{equation}\label{rec}
\Phi(x';q,t)=\Phi(x;q,t) \frac{(-ty;q)_{\infty}}{(y;q)_{\infty}}
\prod_{i=1}^n \frac{(t^2x_iy;q^2)_{\infty}}{(x_iy;q^2)_{\infty}}.
\end{equation}

We will expand both sides of \eqref{rec} in terms of
$P_{\mu}(x;q^2,t^2) y^r$.
After comparing coef\/f\/icients this results in an identity for Pieri
coef\/f\/icients, given in Proposition~\ref{propPieri} below.

\begin{lem}\label{lemma1}
The right-hand side of \eqref{rec} may be expanded as
\[
\sum_{r=0}^{\infty}
\sum_{\substack{\mu,\nu \\ \mu-\nu\in \H}}
\frac{(-t;q)_{r-\abs{\mu-\nu}}}{(q;q)_{r-\abs{\mu-\nu}}}
b^{-}_{\nu}(q,t)
\phi_{\mu/\nu}\big(q^2,t^2\big) P_{\mu}\big(x;q^2,t^2\big) y^r.
\]
\end{lem}

\begin{proof}
By \eqref{Phi}, the $q$-binomial theorem \eqref{qbt}, and the
generating function \eqref{grprod} for the $g_r$, the right
of \eqref{rec} is equal to
\[
\sum_{\nu} \sum_{k,r=0}^{\infty}
\frac{(-t;q)_k}{(q;q)_k}
b^{-}_{\nu}(q,t)
P_{\nu}\big(x;q^2,t^2\big) g_r\big(x;q^2,t^2\big) y^{k+r}.
\]
Recalling the Pieri rule \eqref{P1} this can be further rewritten as
\begin{gather*}
\sum_{k,r=0}^{\infty}   \sum_{\substack{\mu,\nu \\ \mu-\nu\in \H_r}}
\frac{(-t;q)_k}{(q;q)_k}
b^{-}_{\nu}(q,t)
\phi_{\mu/\nu}\big(q^2,t^2\big) P_{\mu}\big(x;q^2,t^2\big) y^{k+r} \\
\qquad{} =\sum_{k=0}^{\infty} \sum_{\substack{\mu,\nu \\ \mu-\nu\in \H}}
\frac{(-t;q)_k}{(q;q)_k}
b^{-}_{\nu}(q,t)
\phi_{\mu/\nu}\big(q^2,t^2\big) P_{\mu}\big(x;q^2,t^2\big) y^{k+\abs{\mu-\nu}} \\
\qquad{} =\sum_{r=0}^{\infty}
\sum_{\substack{\mu,\nu \\ \mu-\nu\in \H}}
\frac{(-t;q)_{r-\abs{\mu-\nu}}}{(q;q)_{r-\abs{\mu-\nu}}}
b^{-}_{\nu}(q,t)
\phi_{\mu/\nu}\big(q^2,t^2\big) P_{\mu}\big(x;q^2,t^2\big) y^r,
\end{gather*}
where the last equality is true since $1/(q;q)_k=0$ for $k$ a negative
integer.
\end{proof}

\begin{lem}\label{lemma2}
The left-hand side of \eqref{rec} may be expanded as
\[
\sum_{r=0}^{\infty}
\sum_{\substack{\la,\mu \\ \la-\mu \in\H_r}}
b^{-}_{\la}(q,t)\psi_{\la/\mu}\big(q^2;t^2\big) P_{\mu}\big(x;q^2,t^2\big)y^r.
\]
\end{lem}

\begin{proof}
Applying \cite[page~348, Example~2]{Macdonald95}
\[
P_{\la}(x';q,t)=\sum_{\substack{\mu \\ \la-\mu\in\H}}
\psi_{\la/\mu}(q,t) P_{\mu}(x;q,t)y^{\abs{\la-\mu}},
\]
the left of \eqref{rec} is equal to
\begin{gather*}
\sum_{\substack{\la,\mu \\ \la-\mu \in\H}}
b^{-}_{\la}(q,t)\psi_{\la/\mu}\big(q^2,t^2\big)
P_{\mu}\big(x;q^2,t^2\big)y^{\abs{\la-\mu}}\\
\qquad{}=\sum_{r=0}^{\infty}\sum_{\substack{\la,\mu \\ \la-\mu \in\H_r}}
b^{-}_{\la}(q,t)\psi_{\la/\mu}\big(q^2,t^2\big) P_{\mu}\big(x;q^2,t^2\big)y^r.
\tag*{\qed}
\end{gather*}\renewcommand{\qed}{}
\end{proof}

Equating the expansions of Lemmas~\ref{lemma1} and \ref{lemma2},
and extracting coef\/f\/icients of $P_{\mu}(x;q^2,t^2)y^r$ it follows
that the proof of Theorem~\ref{thmcon} boils down to a proof of the
following identity for $b^{-}_{\la}(q,t)$.
\begin{prop}\label{propPieri}
For $\mu$ a partition and $r$ a nonnegative integer,
\begin{equation}\label{eqpr}
\sum_{\substack{\nu \\ \mu-\nu\in \H}}
\frac{(-t;q)_{r-\abs{\mu-\nu}}}{(q;q)_{r-\abs{\mu-\nu}}}
b^{-}_{\nu}(q,t)
\phi_{\mu/\nu}\big(q^2,t^2\big)=
\sum_{\substack{\la \\ \la-\mu \in\H_r}}
b^{-}_{\la}(q,t)\psi_{\la/\mu}\big(q^2,t^2\big).
\end{equation}
\end{prop}
Notationally it turns out to be slightly simpler to prove this
in a form involving the Pieri coef\/f\/icient $\psi'_{\la/\mu}(q,t)$
given by
\[
P_{\mu}(x;q,t) e_r(x;q,t)=
\sum_{\substack{\la \\ \la-\mu\in \V_r}} \psi'_{\la/\mu}(q,t) P_{\la}(x;q,t),
\]
where $e_r=P_{(1^r)}$ is the $r$th elementary
symmetric function
\[
e_r(x)=\sum_{i_1<i_2<\dots<i_r} x_{i_1} x_{i_2}\cdots x_{i_r}.
\]
Hence we replace all partitions in \eqref{eqpr} by their conjugates
and use \cite[page~341]{Macdonald95}
\[
\phi_{\mu'/\nu'}(q,t)=\frac{b_{\nu}(t,q)}{b_{\mu}(t,q)}  \psi'_{\mu/\nu}(t,q)
\qquad\text{and}\qquad
\psi_{\la'/\mu'}(q,t)=\psi'_{\la/\mu}(t,q),
\]
as well as \eqref{bp}.
Finally dividing both sides by $b^{-}_{\mu'}(q,t)$
it follows that \eqref{eqpr} can be rewritten as follows.
\addtocounter{thm}{-1}
\begin{prop}\label{prop}\hspace{-3mm}$\mathbf{'}$\;
For $\mu$ a partition and $r$ a nonnegative integer,
\[
\sum_{\substack{\nu \\ \mu-\nu\in \V}}
\frac{(-t;q)_{r-\abs{\mu-\nu}}}{(q;q)_{r-\abs{\mu-\nu}}}
\frac{b^{-}_{\nu}(t,q)}{b^{-}_{\mu}(t,q)}  \psi'_{\mu/\nu}\big(t^2,q^2\big)=
\sum_{\substack{\la \\ \la-\mu \in\V_r}}
\frac{b^{-}_{\la'}(q,t)}{b^{-}_{\mu'}(q,t)}
\psi'_{\la/\mu}\big(t^2,q^2\big).
\]
\end{prop}

Crucial in our proof below is an explicit formula for
$\psi'_{\la/\mu}$ due to Macdonald
\cite[Section~VI.6, Equation~(6.13)]{Macdonald95}
\begin{equation}\label{psipexp}
\psi'_{\la/\mu}(q,t)=
\prod_{\substack{i<j \\ \la_i=\mu_i \\ \la_j=\mu_j+1}}
\frac{(1-q^{\mu_i-\mu_j}t^{j-i-1})(1-q^{\la_i-\la_j}t^{j-i+1})}
{(1-q^{\mu_i-\mu_j}t^{j-i})(1-q^{\la_i-\la_j}t^{j-i})}.
\end{equation}

\begin{proof}[Proof of Proposition \ref{prop}\,$'$]
Let us denote the left- and right-hand sides of the identity
of Proposition~\ref{prop}$'$ by LHS and RHS. Then
\[
\LHS=\sum_{s=0}^r \frac{(-t;q)_s}{(q;q)_s}
\sum_{\substack{\nu \\ \mu-\nu\in \V_{r-s}}}
\frac{b^{-}_{\nu}(t,q)}{b^{-}_{\mu}(t,q)}  \psi'_{\mu/\nu}\big(t^2,q^2\big).
\]
Let $n$ denote the length of the partition $\mu$, i.e.,
$\mu=(\mu_1,\dots,\mu_n)$ with $\mu_n\geq 1$.
We can then replace the sum over $\nu$
by a sum over $k$-subsets $I$
of $[n]$ such that $i\in I$ if\/f $\mu_i-\nu_i=1$.
In other words, $I$ encodes the parts of $\nu$ that dif\/fer from those
of $\mu$. Using this notation as well as \eqref{Iexp} and \eqref{psipexp}
we f\/ind
\begin{gather*}
\frac{b^{-}_{\nu}(t,q)}{b^{-}_{\mu}(t,q)}
=
\prod_{i\in I}
\frac{1-q^{n-i}t^{\mu_i}}{1+q^{n-i+1}t^{\mu_i-1}}
\prod_{\substack{i\in I \\ j\nin I \\ i<j}}
\frac{(1+q^{j-i+1}t^{\mu_i-\mu_j-1})(1-q^{j-i-1}t^{\mu_i-\mu_j})}
{(1-q^{j-i}t^{\mu_i-\mu_j})(1+q^{j-i}t^{\mu_i-\mu_j-1})} \\
\phantom{\frac{b^{-}_{\nu}(t,q)}{b^{-}_{\mu}(t,q)}=}{} \times
\prod_{\substack{i\in I \\ j\nin I \\ i>j}}
\frac{(1+q^{j-i}t^{\mu_i-\mu_j})(1-q^{j-i}t^{\mu_i-\mu_j-1})}
{(1-q^{j-i+1}t^{\mu_i-\mu_j-1})(1+q^{j-i-1}t^{\mu_i-\mu_j})}
\end{gather*}
and
\[
\psi'_{\mu/\nu}(q,t)=\prod_{\substack{i\in I \\ j\nin I \\ i>j}}
\frac{(1-q^{\mu_i-\mu_j-1}t^{j-i+1})(1-q^{\mu_i-\mu_j}t^{j-i-1})}
{(1-q^{\mu_i-\mu_j}t^{j-i})(1-q^{\mu_i-\mu_j-1}t^{j-i})}.
\]
Hence
\begin{gather*}
\sum_{\substack{\nu \\ \mu-\nu\in \V_k}}
\frac{b^{-}_{\nu}(t,q)}{b^{-}_{\mu}(t,q)}  \psi'_{\mu/\nu}\big(t^2,q^2\big) \\
\qquad{}=\sum_{\substack{I\subseteq [n] \\ \abs{I}=k}}
\prod_{i\in I}
\frac{1-q^{n-i}t^{\mu_i}}{1+q^{n-i+1}t^{\mu_i-1}}
\prod_{\substack{i\in I \\ j\nin I}}
\frac{(1+q^{j-i+1}t^{\mu_i-\mu_j-1})(1-q^{j-i-1}t^{\mu_i-\mu_j})}
{(1-q^{j-i}t^{\mu_i-\mu_j})(1+q^{j-i}t^{\mu_i-\mu_j-1})},
\end{gather*}
resulting in
\begin{gather*}
\LHS=\sum_{s=0}^r \frac{(-t;q)_s}{(q;q)_s}
\sum_{\substack{I\subseteq [n] \\ \abs{I}=r-s}}
\prod_{i\in I}
\frac{1-q^{n-i}t^{\mu_i}}{1+q^{n-i+1}t^{\mu_i-1}}
 \prod_{\substack{i\in I \\ j\nin I}}
\frac{(1+q^{j-i+1}t^{\mu_i-\mu_j-1})(1-q^{j-i-1}t^{\mu_i-\mu_j})}
{(1-q^{j-i}t^{\mu_i-\mu_j})(1+q^{j-i}t^{\mu_i-\mu_j-1})}.
\end{gather*}

We next turn to the right-hand side, which is a
sum over partitions $\la$ such that $\la-\mu\in\V_r$.
Recall that $\mu$ has exactly $n$ parts. The maximum number of parts of $\la$
is thus $n+r$ (when $\la=\mu\cup (1^r)$) and the minimum number of parts
of $\la$ is $n$.
Hence we can write
\[
\RHS=
\sum_{s=0}^r \sum_{\substack{\la \\ \la-\mu \in\V_r \\ l(\la)=n+s}}
\frac{b^{-}_{\la'}(q,t)}{b^{-}_{\mu'}(q,t)}
\psi'_{\la/\mu}\big(t^2,q^2\big).
\]
Again we let $I\subseteq [n]$ (with $\abs{I}=r-s$) be the set of indices
of those parts of $\mu$ to which a~square is added to form $\la$; $i\in I$ if\/f
$\la_i-\mu_i=1$ for $i\in[n]$.
For example if $\mu=(3,2,2,1)$ and $\la=(4,3,2,2,1,1)$
then $n=4$, $r=5$, $s=2$ and $I=\{1,2,4\}$.

From \eqref{bprimeexp}, after a tedious calculation, we obtain
\begin{gather*}
\frac{b^{-}_{\la'}(q,t)}{b^{-}_{\mu'}(q,t)}=
\frac{(-t;q)_s}{(q;q)_s}
 \prod_{i\in I} \frac{1+q^{n-i+s}t^{\mu_i+1}}{1-q^{n-i+s+1}t^{\mu_i}}
\prod_{j\nin I} \frac{(1+q^{n-j+s}t^{\mu_j})(1-q^{n-j+1}t^{\mu_j-1})}
{(1-q^{n-j+s+1}t^{\mu_j-1})(1+q^{n-j}t^{\mu_j})}  \\
\phantom{\frac{b^{-}_{\la'}(q,t)}{b^{-}_{\mu'}(q,t)}=}{} \times
\prod_{\substack{i\in I \\ j\nin I \\ i<j}}
\frac{(1+q^{j-i-1}t^{\mu_i-\mu_j+1})(1-q^{j-i+1}t^{\mu_i-\mu_j})}
{(1-q^{j-i}t^{\mu_i-\mu_j})(1+q^{j-i}t^{\mu_i-\mu_j+1})} \\
\phantom{\frac{b^{-}_{\la'}(q,t)}{b^{-}_{\mu'}(q,t)}=}{} \times\prod_{\substack{i\in I \\ j\nin I \\ i>j}}
\frac{(1-q^{j-i}t^{\mu_i-\mu_j+1})(1+q^{j-i}t^{\mu_i-\mu_j})}
{(1-q^{j-i-1}t^{\mu_i-\mu_j+1})(1+q^{j-i+1}t^{\mu_i-\mu_j})}.
\end{gather*}
Furthermore, from \eqref{psipexp},
\begin{gather*}
\psi'_{\la/\mu}(q,t)=
\prod_{j\nin I}
\frac{(1-q^{\mu_j-1}t^{n-j+s+1})(1-q^{\mu_j}t^{n-j})}
{(1-q^{\mu_j}t^{n-j+s})(1-q^{\mu_j-1}t^{n-j+1})} \\
\phantom{\psi'_{\la/\mu}(q,t)=}{} \times
\prod_{\substack{i\in I \\ j\nin I \\ i>j}}
\frac{(1-q^{\mu_i-\mu_j}t^{j-i+1})(1-q^{\mu_i-\mu_j+1}t^{j-i-1})}
{(1-q^{\mu_i-\mu_j+1}t^{j-i})(1-q^{\mu_i-\mu_j}t^{j-i})}.
\end{gather*}
Putting these results together yields
\begin{gather*}
\RHS=\sum_{s=0}^r
\frac{(-t;q)_s}{(q;q)_s}
\sum_{\substack{I\subseteq [n] \\ \abs{I}=r-s}}
\prod_{i\in I} \frac{1+q^{n-i+s}t^{\mu_i+1}}{1-q^{n-i+s+1}t^{\mu_i}}
\prod_{j\nin I} \frac{(1+q^{n-j+s+1}t^{\mu_j-1})(1-q^{n-j}t^{\mu_j})}
{(1-q^{n-j+s}t^{\mu_j})(1+q^{n-j+1}t^{\mu_j-1})} \\
\phantom{\RHS=}{} \times
\prod_{\substack{i\in I \\ j\nin I}}
\frac{(1+q^{j-i-1}t^{\mu_i-\mu_j+1})(1-q^{j-i+1}t^{\mu_i-\mu_j})}
{(1-q^{j-i}t^{\mu_i-\mu_j})(1+q^{j-i}t^{\mu_i-\mu_j+1})}.
\end{gather*}

Finally equating LHS and RHS we obtain
\begin{gather}
\sum_{s=0}^r \frac{(-t;q)_s}{(q;q)_s}
\sum_{\substack{I\subseteq [n] \\ \abs{I}=r-s}}
\prod_{i\in I}
\frac{1-q^{n-i}t^{\mu_i}}{1+q^{n-i+1}t^{\mu_i-1}}
\prod_{\substack{i\in I \\ j\nin I}}
\frac{(1+q^{j-i+1}t^{\mu_i-\mu_j-1})(1-q^{j-i-1}t^{\mu_i-\mu_j})}
{(1-q^{j-i}t^{\mu_i-\mu_j})(1+q^{j-i}t^{\mu_i-\mu_j-1})} \nonumber\\
\qquad{}{} =\sum_{s=0}^r
\frac{(-t;q)_s}{(q;q)_s}
\sum_{\substack{I\subseteq [n] \\ \abs{I}=r-s}}
\prod_{i\in I} \frac{1+q^{n-i+s}t^{\mu_i+1}}{1-q^{n-i+s+1}t^{\mu_i}}
\prod_{j\nin I} \frac{(1+q^{n-j+s+1}t^{\mu_j-1})(1-q^{n-j}t^{\mu_j})}
{(1-q^{n-j+s}t^{\mu_j})(1+q^{n-j+1}t^{\mu_j-1})} \nonumber\\
\qquad\quad{}\times
\prod_{\substack{i\in I \\ j\nin I}}
\frac{(1+q^{j-i-1}t^{\mu_i-\mu_j+1})(1-q^{j-i+1}t^{\mu_i-\mu_j})}
{(1-q^{j-i}t^{\mu_i-\mu_j})(1+q^{j-i}t^{\mu_i-\mu_j+1})}.\label{final}
\end{gather}
This is a limiting case of Theorem~\ref{Thmrn}.
To see this, take the theorem with $p=0$ (recall that $\theta(z;0)=(1-z)$)
and replace the summation index $r$ by $s$. Then carry out the simultaneous
substitutions
\[
(t,v,w,x_i)\mapsto \big(-t,-\epsilon q^{-r}t,q^{-r},q^{i-n}t^{-\mu_i}/\epsilon\big)
\]
(for all $i\in[n]$)
and take the $\epsilon \to 0$ limit.
Finally replacing $s\mapsto r-s$ and multiplying
both sides by $(-t;q)_r/(q;q)_r$ yields \eqref{final}.
\end{proof}

\section{Elliptic hypergeometric series}\label{secell}

\subsection{A new multivariable transformation formula}

To turn the theta-function identity of Theorem~\ref{Thmrn}
into an identity for elliptic hypergeo\-metric series we
apply the well-known procedure of multiple principal
specialisation, see e.g.,~\mbox{\cite{Kajihara04,KN03,Rosengren01,RS05}}.

For $n$, $m$ integers, let
\[
[n]_m:=\{m+1,m+2,\dots,m+n\},
\]
so that $[n]_0=[n]$.
In Theorem~\ref{Thmrn} replace
\begin{gather}
(x_1,x_2,\dots,x_n)
\mapsto \big(t_1,t_1q,\dots,t_1q^{m_1-1},
t_2,t_2q,\dots,t_2q^{m_2-1},
\dots,\nonumber\\
\phantom{(x_1,x_2,\dots,x_n)\mapsto \big(}{}\dots,
t_N,t_Nq,\dots,t_Nq^{m_N-1}\big),\label{MPS}
\end{gather}
where $m_1+\cdots+m_N=n$.
(In the notation of $\lambda$-rings \cite{Lascoux03} we are making the
substitution $x\mapsto \sum\limits_{i=1}^N t_i (1-q^{m_i})/(1-q)$.)

Since $\theta(1)=0$ it follows that
\[
\prod_{\substack{i\in I \\ j\not\in I}} \theta(q x_i/x_j)
\quad\qquad\qquad\qquad
\Biggl(\text{resp.}\quad
\prod_{\substack{i\in I \\ j\nin I}} \theta(q x_j/x_i)\Biggr)
\]
vanishes unless $I$ is of the form
\[
I=\bigcup_{i=1}^N [k_i]_{m_1+\cdots+m_i-k_i}
\qquad\quad
\biggl(\text{resp.}\quad
I=\bigcup_{i=1}^N [k_i]_{m_1+\cdots+m_{i-1}}\biggr),
\]
where $k_1,\dots,k_N$ are integers such that $0\leq k_i\leq m_i$ for each $i$.
Since $\abs{I}=r$ we must of course further impose that
$\abs{k}:=k_1+\cdots+k_N=r$.

The rest is essentially a straightforward calculation, and we only
sketch the details pertaining to the right-hand side of
Theorem~\ref{Thmrn}:
\begin{gather*}
\RHS\xrightarrow{\eqref{MPS}}
\sum_{r=0}^{\abs{m}}   \frac{(v,w)_r}{(qv/t,qw/t)_r}
\Bigl(\frac{q}{t}\Bigr)^{(N+1)r}
\prod_{i=1}^N
\frac{(vt_i/tw,vt_iq^{r-1})_{m_i}}{(vt_i/qw,vt_iq^r/t)_{m_i}}\\
\hphantom{\RHS\xrightarrow{\eqref{MPS}}}{}   \times
\sum_{\substack{k_1,\dots,k_N \\[1pt] \abs{k}=r \\[1pt] 0\leq k_i\leq m_i}}
\prod_{i=1}^N
\frac{(t_i/q,vt_iq^r/t^2,vt_i/qw,vt_iq^r/t)_{k_i}}
{(t_i/t,vt_iq^{r-1}/t,vt_i/tw,vt_iq^{r-1})_{k_i}} \\
\hphantom{\RHS\xrightarrow{\eqref{MPS}}}{}  \times
\prod_{i,j=1}^N
\frac{(q^{-m_j}t_i/t_j,tt_i/t_j)_{k_i}}{(q^{1-m_j}t_i/tt_j,qt_i/t_j)_{k_i}}
\frac{(qt_i/t_j)_{k_i-k_j}}{(tt_i/t_j)_{k_i-k_j}},
\end{gather*}
where we have replaced $n$ by $\abs{m}:=m_1+\cdots+m_N$.
Using
\[
\sum_{r=0}^{\abs{m}}
\sum_{\substack{k_1,\dots,k_N \\[1pt] \abs{k}=r \\[1pt] 0\leq k_i\leq m_i}}
f_{r,k_1,\dots,k_N}=
\sum_{\substack{k_1,\dots,k_N \\[1pt] 0\leq k_i\leq m_i}}
f_{\abs{k},k_1,\dots,k_N}
\]
and making the substitutions
\[
(t,t_i,v,w)\mapsto (b,abqt_i/c,c,d),
\]
we obtain
\[
V_m(a;b,c,d;t)
\prod_{i=1}^N
\frac{(aqt_i/d,abt_i)_{m_i}}{(aqt_i,abt_i/d)_{m_i}},
\]
where, for $t=(t_1,\dots,t_N)$ and $m=(m_1,\dots,m_N)$,
\begin{gather*}
V_m(a;b,c,d;t):=
\sum_{k_1=0}^{m_1}\cdots\sum_{k_N=0}^{m_N}
\prod_{i=1}^N \biggl(\frac{\theta(at_iq^{k_i+\abs{k}})}{\theta(at_i)}
\frac{(abt_i/c,abt_i/d)_{k_i}}{(aqt_i/c,aqt_i/d)_{k_i}} \\
\phantom{V_m(a;b,c,d;t):=}{}  \times
\frac{(at_i,abq^{m_i}t_i)_{\abs{k}}}{(aqt_i/b,aq^{m_i+1}t_i)_{\abs{k}}}
\frac{(aqt_i/b)_{k_i+\abs{k}}}{(abt_i)_{k_i+\abs{k}}}\biggr) \\
\phantom{V_m(a;b,c,d;t):=}{} \times\frac{(c,d)_{\abs{k}}}{(cq/b,dq/b)_{\abs{k}}}
\Bigl(\frac{q}{b}\Bigr)^{(N+1)\abs{k}}
\prod_{i,j=1}^N
\frac{(q^{-m_j}t_i/t_j,bt_i/t_j)_{k_i}}{(q^{1-m_j}t_i/bt_j,qt_i/t_j)_{k_i}}
\frac{(qt_i/t_j)_{k_i-k_j}}{(bt_i/t_j)_{k_i-k_j}}.
\end{gather*}
We note that a particularly succinct way to express this elliptic
hypergeometric series follows by introducing
$k_{N+1}:=-k_1-\cdots-k_N$. Then
\begin{gather*}
V_m(1/t_{N+1};b,c,d;t):=
\sum_{k_1=0}^{m_1}\cdots\sum_{k_N=0}^{m_N}
\prod_{i=1}^{N+1}
\frac{(bt_i/ct_{N+1},bt_i/dt_{N+1})_{k_i}}
{(qt_i/ct_{N+1},qt_i/dt_{N+1})_{k_i}} \\
\phantom{V_m(1/t_{N+1};b,c,d;t):=}{} \times\prod_{i=1}^{N+1}\prod_{j=1}^N
\frac{(q^{-m_j}t_i/t_j,bt_i/t_j)_{k_i}}{(q^{1-m_j}t_i/bt_j,qt_i/t_j)_{k_i}}
\prod_{i,j=1}^{N+1} \frac{(qt_i/t_j)_{k_i-k_j}}{(bt_i/t_j)_{k_i-k_j}}.
\end{gather*}

A similar calculation may be carried out for the left-hand side of
the theorem and we f\/ind exactly the same multiple basic hypergeometric series,
but with $a$ replaced by $\hat{a}:=cd/ab$ and~$t_i$ replaced by
$s_i:=q^{-m_i}/t_i$ for all $i$.
As a result we can claim the following
transformation formula for A$_{N-1}$ elliptic hypergeometric
series.
\begin{thm}\label{thmVst}
Let $\hat{a}=cd/ab$ and $s_it_i=q^{-m_i}$ for all $i\in[N]$.
Then
\[
V_m(a;b,c,d;t)=
V_m(\hat{a};b,c,d;s)
\prod_{i=1}^N \frac{(aqt_i,\hat{a}qs_i/c,\hat{a}qs_i/d,aqt_i/cd)_{m_i}}
{(\hat{a}qs_i,aqt_i/c,aqt_i/d,\hat{a}qs_i/cd)_{m_i}}.
\]
\end{thm}
For $N=1$, after rescaling $(a,\hat{a})\mapsto (a/t_1,\hat{a}/s_1)$
and replacing $m_1\mapsto n$, this gives
\begin{gather}
\sum_{k=0}^n
\frac{\theta(aq^{2k})}{\theta(a)}
\frac{(a,b,c,d,ab/c,ab/d,abq^n,q^{-n})_k}
{(q,aq/b,aq/c,aq/d,cq/b,dq/b,q^{1-n}/b,aq^{n+1})_k}
\frac{(aq/b)_{2k}}{(ab)_{2k}}
\Bigl(\frac{q}{b}\Bigr)^{2k} \nonumber\\
\qquad{}=\frac{(aq,\hat{a}q/c,\hat{a}q/d,aq/cd)_n}
{(\hat{a}q,aq/c,aq/d,\hat{a}q/cd)_n}  \nonumber\\
\quad\qquad{}\times
\sum_{k=0}^n
\frac{\theta(\hat{a}q^{2k})}{\theta(\hat{a})}
\frac{(\hat{a},b,c,d,\hat{a}b/c,\hat{a}b/d,\hat{a}bq^n,q^{-n})_k}
{(q,\hat{a}q/b,\hat{a}q/c,\hat{a}q/d,cq/b,dq/b,q^{1-n}/b,\hat{a}q^{n+1})_k}
\frac{(\hat{a}q/b)_{2k}}{(\hat{a}b)_{2k}}
\Bigl(\frac{q}{b}\Bigr)^{2k},\label{new}
\end{gather}
where $\hat{a}=q^{-n}cd/ab$.
Curiously, even this one-dimensional case, which may also be
written as a transformation between $_{20}V_{19}$ elliptic
hypergeometric series (or for $p=0$ as a transformation between
$_{14}W_{13}$ basic hypergeometric series), is new.
To the best of our knowledge
it is the f\/irst-ever example of a transformation
that does not yield a summation upon specialisation of some of
its parameters. \textbf{We award AU\$25 for a proof}
of \eqref{new} based on known identities for one-variable
elliptic hypergeometric series,
and AU\$10 for a proof of the $p=0$ case using identities
for basic hypergeometric series.
We do remark that \eqref{new} may be viewed as a somewhat strange
generalisation of Jackson's $_6\phi_5$ summation
\cite[Equation~(II.21)]{GR04}.
Indeed, after taking $p=0$ the $b\to 0$ limit can be taken. Close
inspection reveals that on the right the summand vanishes
unless $k=0$, resulting in Jackson's sum
\[
{_6W_5}\big(a;c,d,q^{-n};q,aq^{n+1}/cd\big)=
\frac{(aq,aq/cd;q)_n}{(aq/c,aq/d;q)_n}.
\]
The same in fact applies for general $N$, and setting $p=0$ and
then taking the $b\to 0$ limit in Theorem~\ref{thmVst} leads to the
following A$_{N-1}$ extension of Jackson's $_6\phi_5$ summation.

For $t=(t_1,\dots,t_N)$ let
\[
\Delta(t)=\prod_{1\leq i<j\leq N}(t_i-t_j)
\]
be the Vandermonde product, and set
\[
\frac{\Delta(tq^k)}{\Delta(t)}=
\prod_{1\leq i<j\leq N}
\biggl(\frac{t_iq^{k_i}-t_jq^{k_j}}{t_i-t_j}\biggr).
\]

\begin{cor}\label{cornew}
For $m_1,\dots,m_N$ nonnegative integers, $\abs{m}:=m_1+\cdots+m_N$
and $e_2(k)=\sum\limits_{i<j}k_i k_j$ we have
\begin{gather*}
\sum_{k_1,\dots,k_N\geq 0}
\frac{\Delta(tq^k)}{\Delta(t)}
\prod_{i=1}^N \biggl(\frac{1-at_iq^{k_i+\abs{k}}}{1-at_i}\biggr)
\frac{(c,d;q)_{\abs{k}}}{\prod\limits_{i=1}^N(aqt_i/c,aqt_i/d;q)_{k_i}}
q^{-e_2(k)} \\
 \qquad \quad{} \times \prod_{i=1}^N \Biggl(
\frac{(at_i;q)_{\abs{k}}}{\prod\limits_{j=1}^N (qt_i/t_j;q)_{k_i}}
\frac{\prod\limits_{j=1}^N(q^{-m_j}t_i/t_j;q)_{k_i}}{(aq^{m_i+1}t_i;q)_{\abs{k}}}
\biggl(\frac{at_i q^{\abs{m}+1}}{cd}\biggr)^{k_i} \Biggr)  \\
\qquad{} =\prod_{i=1}^N \frac{(aqt_i,aqt_i/cd;q)_{m_i}}{(aqt_i/c,aqt_i/d;q)_{m_i}}.
\end{gather*}
\end{cor}
This same result also follows by taking the $d\to\infty$ limit in
the $U(n)$ (or $A_{N-1}$) Jackson sum
\cite[Theorem 6.14]{Milne88}, or, alternatively,
by taking the $b\to\infty$ limit in the $D_N$ Jackson sum
\cite[Theorem A.12]{BS98}.

By a standard analytic argument, see e.g., \cite{Rosengren04},
the sum over the $N$-dimensional hyper-rectangle in Theorem~\ref{thmVst}
may be transformed into a sum over the $N$-simplex $k_1,\dots,k_N\geq 0$,
$k_1+\cdots+k_N\leq n$.
\begin{cor}
For $\hat{a}=cdq^{-n}/ab$
\[
W_n(a;b,c,d;s,t)=
W_n(\hat{a};b,c,d;t,s)
\prod_{i=1}^N \frac{(aqt_i,\hat{a}q/dt_i,\hat{a}qs_i/c,aq/cds_i)_n}
{(\hat{a}qs_i,aq/ds_i,aqt_i/c,\hat{a}q/cdt_i)_n},
\]
where
\begin{gather*}
W_n(a;b,c,d;s,t):=
\sum_{\substack{k_1,\dots,k_N\geq 0 \\ \abs{k}\leq n}}
\prod_{i=1}^N \biggl(\frac{\theta(at_iq^{k_i+\abs{k}})}{\theta(at_i)}
\frac{(abt_i/c,abq^nt_i)_{k_i}}{(aqt_i/c,aq^{n+1}t_i)_{k_i}} \\
\phantom{W_n(a;b,c,d;s,t):=}{}  \times
\frac{(at_i,ab/ds_i)_{\abs{k}}}{(aqt_i/b,aq/ds_i)_{\abs{k}}}
\frac{(aqt_i/b)_{k_i+\abs{k}}}{(abt_i)_{k_i+\abs{k}}}\biggr) \\
\phantom{W_n(a;b,c,d;s,t):=}{} \times\frac{(c,q^{-n})_{\abs{k}}}{(cq/b,q^{1-n}/b)_{\abs{k}}}
\Bigl(\frac{q}{b}\Bigr)^{(N+1)\abs{k}}
\prod_{i,j=1}^N
\frac{(dt_is_j,bt_i/t_j)_{k_i}}{(dqt_is_j/b,qt_i/t_j)_{k_i}}
\frac{(qt_i/t_j)_{k_i-k_j}}{(bt_i/t_j)_{k_i-k_j}}.
\end{gather*}
\end{cor}

{}From the $p=0$ case of Theorem~\ref{thmVst} we may also deduce the following
double multi-sum identity.
\begin{cor}\label{cordmsum}
For $m_1,\dots,m_N$ nonnegative integers and $\abs{m}=m_1+\dots+m_N$ we have
\begin{gather*}
 \sum_{\substack{l_i,k_i\geq 0\\[1pt]l_i+k_i\le m_i\\[1pt] i\in[N]}}
\prod_{i=1}^N \biggl(\frac{(at_iq^{\abs{m}};q)_{l_i+k_i}}
{(abt_iq^{\abs{m}};q)_{l_i+k_i}}
\frac{(aqt_i/cd;q)_{l_i+k_i-\abs{k}}}{(abt_i/cd;q)_{l_i+k_i-\abs{k}}}
\frac{(abt_i/c,abt_i/d;q)_{l_i}}{(aqt_i/c,aqt_i/d;q)_{l_i}}\\
\qquad\quad{} \times
\frac{(abt_i-cdq^{\abs{k}-l_i})}
{(abt_i-cdq^{\abs{k}-l_i-k_i})}
\frac{(cd/at_i;q)_{\abs{k}}}{(cdq/abt_i;q)_{\abs{k}}}
\frac{(cdq^{1-l_i+\abs{k}-k_i}/ab^2t_i;q)_{k_i}}
{(cdq^{-l_i+\abs{k}-k_i}/at_i;q)_{k_i}}\biggr)\\
\qquad\quad{} \times
\prod_{i,j=1}^N\frac{(q^{-m_j}t_i/t_j;q)_{l_i+k_i}}
{(q^{-m_j}bt_i/t_j;q)_{l_i+k_i}} \frac{(bt_i/t_j;q)_{l_i}}
{(qt_i/t_j;q)_{l_i}} \frac{(q^{l_i-l_j}bt_i/t_j;q)_{k_i}}
{(q^{1+l_i-l_j}t_i/t_j;q)_{k_i}}\\
\qquad\quad{} \times
\frac{(c,d;q)_{\abs{k}}}{(cq/b,dq/b;q)_{\abs{k}}}
q^{\abs{l}+N\abs{k}}b^{(1-N)\abs{k}}
\frac{1}{\Delta(t)} \\
\qquad\quad{} \times
\det_{1\le i,j\le N}\!\bigg[(t_iq^{l_i+k_i})^{N-j}\bigg(\!1-
b^{N-j+1}\frac{(1-at_iq^{l_i+k_i+\abs{m}})}{(1-abt_iq^{l_i+k_i+\abs{m}})}
\prod_{r=1}^N\frac{(t_iq^{l_i+k_i}-t_rq^{m_r})}
{(bt_iq^{l_i+k_i}-t_rq^{m_r})}\bigg)\bigg]\\
\qquad =q^{\sum\limits_{i=1}^N(i+1)m_i} b^{-2\abs{m}}
\prod_{i,j=1}^N\frac{(bt_i/t_j;q)_{m_i}}{(qt_i/bt_j;q)_{m_i}}
\prod_{1\leq i<j<N}\frac{(qt_i/bt_j;q)_{m_i-m_j}}{(bt_i/t_j;q)_{m_i-m_j}} \\
\qquad\quad{} \times
\frac{(c,d;q)_{\abs{m}}}{(cq/b,dq/b;q)_{\abs{m}}}
\prod_{i=1}^N\frac{(abt_i/c,abt_i/d;q)_{m_i}}{(aqt_i/c,aqt_i/d;q)_{m_i}}
\frac{(abt_i;q)_{\abs{m}}}{(aqt_i/b;q)_{\abs{m}}}
\frac{(aqt_i/b;q)_{m_i+\abs{m}}}{(abt_i;q)_{m_i+\abs{m}}}.
\end{gather*}
\end{cor}

\begin{proof}
We def\/ine
\begin{gather*}
f_k(a;b,c,d;t)=
q^{N|k|}b^{-(N+1)\abs{k}}
\prod_{i=1}^N\frac{(abt_i/c,abt_i/d;q)_{k_i}}{(aqt_i/c,aqt_i/d;q)_{k_i}}
\frac{(at_i;q)_{\abs{k}}}{(aqt_i/b;q)_{\abs{k}}}
\frac{(aqt_i/b;q)_{k_i+\abs{k}}}{(abt_i;q)_{k_i+\abs{k}}}\\
\phantom{f_k(a;b,c,d;t)=}{} \times
\frac{(c,d;q)_{\abs{k}}}{(cq/b,dq/b;q)_{\abs{k}}}
\prod_{i,j=1}^N\frac{(bt_i/t_j;q)_{k_i}}{(qt_i/t_j;q)_{k_i}}
\frac{(qt_i/t_j;q)_{k_i-k_j}}{(bt_i/t_j;q)_{k_i-k_j}}
\end{gather*}
and
\begin{gather*}
M_{mk}(a;b,c,d;t)=q^{\abs{k}}
\prod_{i=1}^N\frac{(1-at_iq^{k_i+\abs{k}})}{(1-at_i)}
\frac{(abq^{m_i}t_i;q)_{\abs{k}}}{(aq^{1+m_i}t_i;q)_{\abs{k}}}
\frac{(aqt_i/c,aqt_i/d;q)_{m_i}}{(aqt_i,aqt_i/cd;q)_{m_i}}\\
\phantom{M_{mk}(a;b,c,d;t)=}{} \times
\prod_{i,j=1}^N\frac{(q^{-m_j}t_i/t_j;q)_{k_i}}
{(q^{1-m_j}t_i/bt_j;q)_{k_i}}.
\end{gather*}
Then the $p=0$ case of Theorem \ref{thmVst} may be put as
\begin{equation*}
\sum_{\substack{0\le k_i\le m_i\\[1pt]i\in[N]}}
M_{mk}(a;b,c,d;t) f_k(a;b,c,d;t)=
\sum_{\substack{0\le k_i\le m_i\\[1pt]i\in[N]}}
M_{mk}(\hat{a};b,c,d;s) f_k(\hat{a};b,c,d;s),
\end{equation*}
where $\hat{a}=cd/ab$ and $s_it_i=q^{-m_i}$ for all $i\in[N]$.
According to \cite[Theorem 2.3]{LS06} (modulo the simultaneous
substitutions $b\mapsto t_1\cdots t_N/a$, $a_i(y_i)\mapsto bt_iq^{y_i}$,
$c_i(y_i)\mapsto t_i q^{y_i}$ for $i\in{[N]}$)
the inverse of the inf\/inite-dimensional lower-triangular matrix $M_{mk}$
is given by
\begin{gather*}
M_{mk}^{-1}(a;b,c,d;t)=
\frac{\Delta(tq^m)}{\Delta(t)^2}
\prod_{i=1}^N\frac{(at_i;q)_{k_i+\abs{m}}}{(abt_i;q)_{k_i+\abs{m}}}
\frac{(abt_i,aqt_i/cd;q)_{k_i}}{(aqt_i/c,aqt_i/d;q)_{k_i}}
\\ \qquad\quad{}  \times
q^{\abs{k}-\abs{m}}
\prod_{i,j=1}^N\frac{(q^{-m_j}t_i/t_j,bt_i/t_j;q)_{k_i}}
{(qt_i/t_j,q^{-m_j}bt_i/t_j;q)_{k_i}}
\frac{(qt_i/bt_j;q)_{m_i}}{(qt_i/t_j;q)_{m_i}} \\
\qquad\quad{}\times
\det_{1\le i,j\le N}\!\bigg[(t_iq^{k_i})^{N-j}\bigg(1-
b^{N-j+1}\frac{(1-at_iq^{k_i+\abs{m}})}{(1-abt_iq^{k_i+\abs{m}})}
\prod_{r=1}^N\frac{(t_iq^{k_i}-t_rq^{m_r})}
{(bt_iq^{k_i}-t_rq^{m_r})}\bigg)\bigg].
\end{gather*}
Hence the $p=0$ case of Theorem \ref{thmVst} is equivalent to
\begin{equation*}
\sum_{\substack{0\le k_i\le m_i\\[1pt] i\in[N]}}
\sum_{\substack{k_i\le l_i\le m_i\\[1pt] i\in[N]}}
M_{ml}^{-1}(a;b,c,d;t)
M_{lk}(\hat{a};b,c,d;s) f_k(\hat{a};b,c,d;s)=f_m(a;b,c,d;t),
\end{equation*}
where $\hat{a}=cd/ab$ and $s_it_i=q^{-l_i}$ for all $i\in[N]$.
Using
\[
\sum_{\substack{0\le k_i\le m_i\\[1pt]k_i\le l_i\le m_i\\[1pt]i\in[N]}}
c_{m,l,k}=
\sum_{\substack{l_i,k_i\ge 0\\[1pt]l_i+k_i\le m_i\\[1pt]i\in[N]}}
c_{m,l+k,k}
\]
and the explicit forms for $f_k$, $M_{mk}$ and $M_{mk}^{-1}$ this
gives (after some elementary manipulations) the result as stated
in the corollary.
\end{proof}

We remark that for $N=1$ the identity in Corollary \ref{cordmsum}
(in which case the determinant appearing in the summand of the
double sum factorises)
admits the following elliptic extension:
\begin{gather*}%\label{doublesum}\notag
\sum_{l,k\geq 0}
\frac{\theta(abq^{2l+2k})}{\theta(ab)}
\frac{(aq^m,q^{-m})_{l+k}}{(bq^{1-m},abq^{m+1})_{l+k}}
\frac{(b,ab/c,ab/d,aq/cd)_l}{(q,aq/c,aq/d,ab/cd)_l}  q^l \\
\qquad \quad{} \times \frac{\theta(cdq^{k-l}/ab)}{\theta(cdq^{-l}/ab)}
\frac{(cdq^{1-l}/ab^2)_k}{(cdq^{-l}/a)_k}
\frac{(b,c,d,cd/a)_k}{(q,cq/b,dq/b,cdq/ab)_k}   q^k \\
\qquad{} =\frac{(aq/b)_{2m}}{(ab)_{2m}}
\frac{(abq,b,c,d,ab/c,ab/d)_m}
{(1/b,aq/b,aq/c,aq/d,cq/b,dq/b)_m}
\left(\frac{q}{b^2}\right)^m.
\end{gather*}
This can proved by inverting the transformation in \eqref{new}
using the elliptic
matrix inversion of \cite[Equation~(3.5); $r=1$]{Warnaar02}.
Since the explicit form of the elliptic extension of the
multivariate matrix inversion of \cite[Theorem~2.3]{LS06}
is not (yet) available, we were unable to extend our basic
($p=0$) double multi-sum identity in Corollary~\ref{cordmsum}
to the elliptic case.

\subsection[Relation to the $q,t$-Littlewood-Richardson coefficients]{Relation to the $\boldsymbol{q,t}$-Littlewood--Richardson coef\/f\/icients}

The $p=0$ instance of the transformation in Theorem \ref{thmVst}
is closely related to an identity that arises by exploiting
some basic symmetries of the $q,t$-Littlewood--Richardson (LR)
coef\/f\/icients. The latter are def\/ined by \cite[Section~VI.7]{Macdonald95}
\begin{equation}\label{PPP}
P_{\mu}(x;q,t)P_{\nu}(x;q,t)=
\sum_{\la} f_{\mu\nu}^{\la}(q,t) P_{\la}(x;q,t).
\end{equation}
Iteration of \eqref{PPP} yields
\begin{equation*}
P_{\mu}(x;q,t)P_{\nu}(x;q,t)P_{\rho}(x;q,t)=
\sum_{\la,\tau} f_{\mu\nu}^{\la}(q,t)
f_{\la\rho}^{\tau}(q,t) P_{\tau}(x;q,t).
\end{equation*}
Since the left-hand side is symmetric under interchange of $\nu$ and $\rho$,
we also have
\begin{equation*}
P_{\mu}(x;q,t)P_{\nu}(x;q,t)P_{\rho}(x;q,t)=
\sum_{\la,\tau} f_{\mu\rho}^{\la}(q,t)
f_{\la\nu}^{\tau}(q,t) P_{\tau}(x;q,t).
\end{equation*}
By equating coef\/f\/icients of $P_{\tau}(x;q,t)$
we immediately get the general identity
\begin{equation}\label{ffff}
\sum_{\la} f_{\mu\nu}^{\la}(q,t) f_{\la\rho}^{\tau}(q,t)=
\sum_{\la} f_{\mu\rho}^{\la}(q,t) f_{\la\nu}^{\tau}(q,t).
\end{equation}
We now assume that $\nu$ and $\rho$ are partitions containing
a single part only. Then the four \mbox{$q,t$-LR} coef\/f\/icients in the above
transformation all become Pieri coef\/f\/icients, and hence completely
factorise. Specif\/ically we set
$\mu=(\mu_1,\dots,\mu_N)$, $\nu=(r)$ and $\rho=(s)$. Using
\eqref{P1} the transformation \eqref{ffff} then simplif\/ies to
\begin{equation}\label{pppp}
\sum_{\substack{\la\\[1pt]\tau-\la\in\H_s\\[1pt]\la-\mu\in\H_r}}
\phi_{\tau/\la}(q,t)\phi_{\la/\mu}(q,t)=
\sum_{\substack{\la\\[1pt]\tau-\la\in\H_r\\[1pt]\la-\mu\in\H_s}}
\phi_{\tau/\la}(q,t)\phi_{\la/\mu}(q,t).
\end{equation}
If we f\/inally apply Macdonald's formula
\cite[Section~VI.6, Example~2(a)]{Macdonald95}
\begin{equation*}
\phi_{\la/\mu}(q,t)=\prod_{1\le i<j\le l(\la)}
\frac{f(q^{\la_i-\la_j}t^{j-i})f(q^{\mu_i-\mu_{j+1}}t^{j-i})}
{f(q^{\la_i-\mu_j}t^{j-i})f(q^{\mu_i-\la_{j+1}}t^{j-i})},
\end{equation*}
where $f(u):=(tu;q)_\infty/(qu;q)_\infty$, we obtain a transformation
for multiple basic hypergeometric series.
To make this explicit we note that since $l(\mu)\leq N$ we have
$l(\la)\leq N+1$ and $l(\tau)\leq N+2$.
We also have $\abs{\tau}-\abs{\mu}=s+r$ so that $\tau_{N+2}$ may
be eliminated by
$\tau_{N+2}=\abs{\mu}+r+s-(\tau_1+\dots+\tau_{N+1})$.
Finally we note that since $\abs{\la}=\abs{\mu}+r$ on the left and
$\abs{\la}=\abs{\mu}+s$ on the right, we may in a similar manner
eliminate the summation index $\la_{N+1}$ on both sides.
If we then shift the remaining summation indices $\la_1,\dots,\la_N$
by $\la_i\mapsto \la_i+\mu_i$ (for $i\in[N]$) and carry out some
simplif\/ications we arrive at the $p=0$ case of
Theorem~\ref{thmVst} subject to the simultaneous substitutions
\begin{equation*}
(a,b,c,d,t_i,m_i)\mapsto
\big(q^{-r},t,q^{\tau_{N+1}-r}t,q^{\abs{\mu}+s-(\tau_1+\dots+\tau_{N+1})},
q^{\mu_i}t^{N+1-i},\tau_i-\mu_i\big)
\end{equation*}
for $i\in[N]$.

We further remark that the identity in
Corollary~\ref{cordmsum} is related to a
recursion formula for the $q,t$-Littlewood--Richardson
coef\/f\/icients given in \cite{Schlosser06}. More precisely,
\cite[Theorem~2.1]{Schlosser06} describes an explicit recursion
in $n$ for the $q,t$-LR coef\/f\/icient
$f_{\mu'\nu'}^{\la'}(t,q)$ where $\la$, $\mu$, and $\nu$
are three partitions with $l(\mu)\leq m$, $l(\nu)\leq n+1$, and
$\abs{\mu}+\abs{\nu}=\abs{\la}$. Choosing $m=1$ this
coef\/f\/icient reduces to a Pieri coef\/f\/icient and hence factorises.
As a result one obtains an identity
equivalent to Corollary \ref{cordmsum}.
This should come as no surprise since \cite[Theorem~2.1]{Schlosser06}
is derived in exactly the same way as Corollary \ref{cordmsum} above;
by the use of multidimensional inverse relations.

\subsection*{Acknowledgements}
We thank the anonymous referees for helpful comments.
The work reported in this paper is supported by the
Australian Research Council and by FWF Austrian Science Fund grants
\hbox{P17563-N13} and S9607.

\pdfbookmark[1]{References}{ref}
\LastPageEnding

\end{document}